\documentclass[11pt]{amsart}
\usepackage{amssymb,latexsym}
\newtheorem{theorem}{Theorem}[section]
\newtheorem{lemma}[theorem]{Lemma}
\newtheorem{proposition}[theorem]{Proposition}

\theoremstyle{definition}

\theoremstyle{remark}
\newtheorem{remark}[theorem]{Remark}

\numberwithin{equation}{section}
\begin{document}
\title[semilinear problem on Heisenberg group]{ Liouville theorem  for a class semilinear elliptic problem on Heisenberg group}
\author{Xi-Nan Ma}
\address{Department of Mathematics\\
         University of Science and Technology of China\\
         Hefei, 230026, Anhui Province, China.}
\email{xinan@ustc.edu.cn}
\thanks{
The research of the first author was supported by NSFC  11871255 and NSFC 11721101.
Research of the second author was supported by NSFC 11861016.\\}
%\footnote{Corresponding author: Qianzhong Ou.}
\author{Qianzhong Ou}
\address{ School of Mathematics and Statistics\\
          Guangxi Normal University\\
           Guilin, 541004, Guangxi Province, China.}
\email{ouqzh@gxnu.edu.cn}
 \maketitle
\begin{abstract}
We obtain an entire Liouville type theorem  to the classical semilinear  subcritical elliptic equation on Heisenberg group.
A pointwise estimate near the isolated singularity was also proved.
The soul of the proofs is an {\it a priori} integral estimate,
which deduced from a generalized formula of that found by Jerison and Lee.
\end{abstract}

%{\bfseries Key words}\quad Heisenberg group, CR Yamabe problem, Liouville theorem

 %\large\normalsize

\section{Introduction}

\setcounter{equation}{0}

In this paper, we study the following equation

\begin{equation}\label{1.1}
   -\triangle_{\mathbb{H}^n}u =2n^2u^{q} \quad \text{in} \quad \Omega,
\end{equation}
where $\Omega$ is a domain in Heisenberg group $\mathbb{H}^n$, and  $u$ is a smooth, nonnegative real function defined in $\Omega$,
while $\triangle_{\mathbb{H}^n}u=u_{\alpha\overline{\alpha}}+u_{\overline{\alpha}\alpha}$ is the Heisenberg laplacian of $u$.
Let $Q=2n+2$ be the homogeneous dimension of  $\mathbb{H}^n$. Denote $q_*=\frac{Q}{Q-2}$ and $q^*=\frac{Q+2}{Q-2}$.  We will deduce
an entire Liouville type theorem  and a point wise estimate near the isolated singularity for  solutions to (\ref{1.1}).  Precisely, we have

\begin{theorem}\label{Thm1}
Let $\Omega=\mathbb{H}^n$ be the whole space and $1<q<q^*$, then the equation (\ref{1.1}) has no positive  solution,
namely, any nonnegative entire solution of  (\ref{1.1}) must be the trivial one.
\end{theorem}

\begin{theorem}\label{Thm2}
Let $\Omega=B_1(0)\backslash\{0\}$ be the punctured unit ball in $\mathbb{H}^n$  and $1<q<q^*$,
then any positive solution $u$ of (\ref{1.1}) satisfies:

\begin{equation}\label{1.2}
   u(\xi)\leq C|\xi|^{\frac{-2}{q-1}} \quad \text{for} \,\, \xi\,\, \text{near}\,\, 0,
\end{equation}
with some positive constant $C$ depending only on $n$ and $q$.
\end{theorem}

The soul of the proofs of theorem \ref{Thm1},\ref{Thm2} is an integral estimate, which may be interested itself.
In fact we shall prove the following

\begin{theorem}\label{Thm3}
 Let $1<q<q^*$, $B_{4r}(\xi_0)\subset \Omega$ be any ball centered at $\xi_0$ with radio $4r$.
 Then any positive solution $u$ of (\ref{1.1}) satisfies:

\begin{equation}\label{1.3}
  \int_{B_r(\xi_0)}  u^{3q-q^*} \leq C\,r^{Q-2\times\frac{3q-q^*}{q-1}} ,
\end{equation}
with some positive constant $C$ depending only on $n$ and $q$.
\end{theorem}

For  $1<q<q^*$, we see $Q-2\times\frac{3q-q^*}{q-1}<0$. So if $u$ be a positive solution of (\ref{1.1}) with $\Omega=\mathbb{H}^n$,
taking  $r\rightarrow +\infty$ in  (\ref{1.3}) we have

\begin{equation}\label{1.4}
  \int_{\mathbb{H}^n}  u^{3q-q^*} \leq 0.
\end{equation}
This contradiction signifies directly the conclusion of theorem \ref{Thm1}.
Also, we will prove theorem \ref{Thm2} by using (\ref{1.3})
combining with the Harnack inequality deduced by Capogna-Danielli-Garofalo ( see Theorem 3.1 in \cite{CDG1993}).

The equation (\ref{1.1}) had been studied intensively by many authors in decades. In fact, it comes from the CR Yamabe problem on $\mathbb{H}^n$.
Let $\mathbf{\Theta}$ be the standard contact form on $\mathbb{H}^n$. Consider another smooth contact form $\theta=u^{\frac{2}{n}}\mathbf{\Theta}$.
Then the pseudo-Hermitian scalar curvature associated to the Fefferman metric of ($\mathbb{H}^n$,$\theta$) is $R=4n(n+1)u^{q-q^*}$ while $u$ satisfies the  equation (\ref{1.1}). Especially, for $q=q^*$, the pseudo-Hermitian scalar curvature $R$ is a constant and it is called  CR Yamabe problem to find such a contact form $\theta$. Accordingly, for $q=q^*$, the equation (\ref{1.1}) is called the CR Yamabe equation. The constant $1+q^*=\frac{2Q}{Q-2}$ is the CR Sobolev embedding exponent and, for the Yamabe equation, there is nontrivial solution as follows

\begin{equation}\label{1.5}
 u(z,t)=C\big|t+\sqrt{-1}z\cdot \overline{z}+z\cdot \mu +\lambda\big|^{-n}
\end{equation}
for some $C>0$, $\lambda\in \mathbf{C}$, Im($\lambda$)$>|\mu|^2/4$, and $\mu\in \mathbf{C}^n$, which is also the only extremals of
the CR Yamabe functional, or the CR Sobolev inequality, on $\mathbb{H}^n$.  So our theorem 1 also confirms that $q^*$ is really  critical.
The equation  (\ref{1.1}) had catched  many mathematician's attention since it raised in the CR Yamabe problem. The CR Yamabe problem had been initiated and studied by David Jerison and John M. Lee in their series fundamental works (see \cite{JL1987}-\cite{JL1989}).
For compact, strictly pseudovonvex CR manifold, the CR Yamabe problem had been solved in case of not locally CR equivalent to sphere $\mathbf{S}^{2n+1}$  by Jerison-Lee \cite{JL1989} for $n\geq 2 $ and Gamara \cite{Ga2001} for $n=1$, and in case of locally
CR equivalent to $\mathbf{S}^{2n+1}$ by Gamara-Yacoub \cite{GY2001} for all $n\geq 1$. The CR Yamabe problem on closed Einstein pseudohermitian manifold was also studied by Wang \cite{Wang2015}.

On $\mathbb{H}^n$, the uniqueness of CR Yamabe solutions  was also obtained by Jerison-Lee \cite{JL1988} for the case of finite volume, i.e., $u\in L^{\frac{2Q}{Q-2}}(\mathbb{H}^n)$, and by Garofalo-Vassilev \cite{GV2001} for the case of cylindrically symmetry on groups of Heisenberg type.  For the subcritical case $1<q\leq q_*$, Birindelli-Dolcetta \cite{BDC1997} proved that the only nonnegative entire solution of (\ref{1.1}) is the trivial one, where they also showed that $q=q_*$ is sharp for the nonexistence of the inequality

\begin{equation}\label{1.6}
   -\triangle_{\mathbb{H}^n}u \geq 2n^2u^{q} \quad \text{in} \quad \mathbb{H}^n.
\end{equation}
For the subcritical case $q_*<q<q^*$, the classification of solutions to the equations (\ref{1.1}) is still open, except
for some partial results, such as the solutions are cylindrical or decay at infinity in \cite{BP1999},
and as $n>1,\,1<q\leq  q^*-\frac{1}{(Q-2)(Q-1)^2}$ in \cite{Xu2009}.

  There are analogous results in the Euclidean case. In the splendid paper \cite{GS1981}, B. Gidas and J. Spruck proved that,
for $1<q<\frac{n+2}{n-2}$, the following equation (\ref{1.7}) has no positive entire solution
in the n-dimension Euclidean space $\mathbb{R}^n$:

\begin{equation}\label{1.7}
   -\triangle u = u^{q} .
\end{equation}
The method used by Gidas-Spruck \cite{GS1981} is the integral estimate, as here in our paper.
Later, Chen-Li [3] also got the same result by using the method of moving plane. Gidas-Spruck \cite{GS1981}
also gave a singularity estimate, precisely, they proved that for $\frac{n}{n-2}<q<\frac{n+2}{n-2}$,
the positive solution of (\ref{1.7}) in the punctured unit ball, with a nonremovable  singularity at the origin, must satisfies

\begin{equation}\label{1.8}
   |x|^{\frac{2}{q-1}}u(x)\rightarrow C_0 \quad \text{as} \,\, x\rightarrow 0.
\end{equation}
Also as in the Euclidean case, the Liouville type result in theorem \ref{Thm1} may be useful
in resolving the Dirichlet problem of the same equation, via the blow-up analysis, .

To get the integral estimate (\ref{1.3}), there are usually two difficulties to be overcome in a noncompact domain.
One is to find a suitable identity, and the other is to estimate the ``tail" terms after integrating
by part of the identity multiplied  suitable cut-off function. When they studied the CR Yamabe problem in their splendid
work \cite{JL1988}, Jerison and Lee  had found several remarkable identities  with the help of a computer program.
 The idea of Jerison-Lee \cite{JL1988} was originally due to Obata in his classic work \cite{Ob1971}. Roughly speaking,
 the main idea is to find an identity to express some suitable nonnegative terms (usually with associate geometry data)
 in a divergence form. Then integrating both sides of the identity to get useful results as one desired.
 But there is a pity that the identity in case of Heisenberg group $\mathbb{H}^n$
 given by  Jerison-Lee(see (4.2) for example in \cite{JL1988})  is in a so complicated form that,
one must suffer an awful long and tedious computation to check it and, may hardly to make any generalization.
Nevertheless, based on our new observation, we would generalize the Jerison-Lee's identity to a new form
and gave a transparent proof, so that it can be used in dealing with the subcritical case of the equation (\ref{1.1}).

The paper is organized as follows. In section 2, we introduced  some  notations and proved a generalization of the Jerison-Lee's identity.
Then, using this generalized identity, we proved theorem 3 in section 3. The proof of theorem 2 shall be presented in section 4.

\section{Generalization of Jerison-Lee's identity }

\setcounter{equation}{0}
\setcounter{theorem}{0}

In this Section we discuss the generalization of a remarkable Jerison-Lee's identity from \cite{JL1988} on Heissenberg group  $\mathbb{H}^{n}$.
We adopt notations as in \cite{JL1988}.

We shall first give a brief introduction to the Heissenberg group  $\mathbb{H}^{n}$ and some notations.
   We consider  $\mathbb{H}^{n}$  as the set $\mathbb{C}^n\times \mathbb{R}$ with coordinates ($z,\,t$) and group law $\circ$:

   $$(z,t)\circ(\xi,t)=\big( z+\xi,\, t+s+2\mathbf{Im} z^{\alpha}\overline{\xi}^{\alpha}\big)\quad
     \text{for}\,\,(z,t),\,(\xi,t)\in \mathbb{C}^n\times \mathbb{R}, $$
\noindent where and in the sequel, the repeated indices are sum form $1$ to $n$. The CR structure of $\mathbb{H}^{n}$ is given by the bundle $\mathcal{H}$ spanned by the left-invariant vector
fields $Z_{\alpha} = \partial/\partial z^{\alpha}+ \sqrt{-1}\overline{z}^{\alpha} \partial/\partial t$, $\alpha= 1, \cdots, n$.
The standard (left-invariant) contact form on $\mathbb{H}^{n}$ is
$\mathbf{\Theta}= dt + \sqrt{-1}(z^{\alpha}d\overline{z}^{\alpha} - \overline{z}^{\alpha}dz^{\alpha})$.
With respect to the standard holomorphic frame $\{Z_{\alpha}\}$ and dual admissible coframe $\{dz^{\alpha}\}$,
the  $Levi\,\, forms\,\,h_{\alpha\overline{\beta}}= 2\delta_{\alpha\overline{\beta}} $.
Accordingly, for a smooth function $f$  on $\mathbb{H}^{n}$, denote its derivatives by
$f_{\alpha}= Z_{\alpha}f$, $f_{\alpha\overline{\beta}}= Z_{\overline{\beta}}(Z_{\alpha}f)$
, $f_0= \frac{\partial f}{\partial t}$, $f_{0\alpha}= Z_{\alpha}(\frac{\partial f}{\partial t})$, etc.
We would also indicate the derivatives of functions or vector fields with indices preceded by a comma, to avoid confusion.
Then we have the following commutative formulae:

$$ f_{\alpha\beta}-f_{\beta\alpha}=0,\quad
 f_{\alpha\overline{\beta}}-f_{\overline{\beta}\alpha}=2\sqrt{-1}\delta_{\alpha\overline{\beta}}\,f_0,\quad
 f_{0\alpha}-f_{\alpha0}=0,$$
$$f_{\alpha\beta 0}-f_{\alpha 0\beta}=0,\qquad
f_{\alpha\beta\overline{\gamma}}-f_{\alpha\overline{\gamma}\beta}
=2\sqrt{-1}\delta_{\beta\overline{\gamma}}\,f_{\alpha 0},\, \cdots.$$

Now we are at the point to give the  generalized identity for positive solution of the equation (\ref{1.1}).
Let $u>0$ solves (\ref{1.1}). Take $e^f= u^{\frac{1}{n}}$ and $q=q^*+ \frac{p}{n}$, then $f$ satisfies the following equation

\begin{equation}\label{2.1}
\mathbf{Re} f_{\alpha\overline{\alpha}}=-n|\partial f|^2-ne^{(2+p)f},
\end{equation}
where $|\partial f|^2=f_{\alpha}f_{\overline{\alpha}}$.
Define the tensors

\begin{equation}\label{2.2}
\begin{split}
 D_{\alpha\beta}            =& f_{\alpha\beta}-2f_{\alpha}f_{\beta}, \qquad\qquad\qquad  D_{\alpha}=D_{\alpha\beta}f_{\overline{\beta}},\\
 E_{\alpha\overline{\beta}} =& f_{\alpha\overline{\beta}}-\frac{1}{n}f_{\gamma\overline{\gamma}}\delta_{\alpha\overline{\beta}},
                              \qquad\quad \qquad E_{\alpha}=E_{\alpha\overline{\beta}}f_{\beta}, \\
 G_{\alpha}=& \sqrt{-1}f_{0\alpha}-\sqrt{-1}f_0f_{\alpha}+e^{(2+p)f}f_{\alpha}+|\partial f|^2f_{\alpha}.
\end{split}
\end{equation}

Denote the function $g=|\partial f|^2+e^{(2+p)f}-\sqrt{-1}f_0$. Then we can rewrite the equation (\ref{2.1}) as

\begin{equation}\label{2.3}
 f_{\alpha\overline{\alpha}}=-ng.
\end{equation}

\noindent Moreover, we  observe that

\begin{equation}\label{2.4}
\begin{split}
\,&  E_{\alpha\overline{\beta}} = f_{\alpha\overline{\beta}}+g\delta_{\alpha\overline{\beta}},\qquad\,\,\qquad
     E_{\alpha}= f_{\alpha\overline{\beta}}f_{\beta}+gf_{\alpha},\\
\,&  D_{\alpha}=f_{\alpha\beta}f_{\overline{\beta}}-2|\partial f|^2f_{\alpha},\qquad
     G_{\alpha}=\sqrt{-1}f_{0\alpha}+gf_{\alpha},
\end{split}
\end{equation}
and by

\begin{equation}\label{2.5}
 (|\partial f|^2)_{,\overline{\alpha}} =D_{\overline{\alpha}}+E_{\overline{\alpha}}+\overline{g}f_{\overline{\alpha}}-2f_{\overline{\alpha}}e^{(2+p)f},
\end{equation}
we find

\begin{equation}\label{2.6}
\begin{split}
g_{\overline{\alpha}} =& D_{\overline{\alpha}}+E_{\overline{\alpha}}+G_{\overline{\alpha}}+pf_{\overline{\alpha}}e^{(2+p)f}.
\end{split}
\end{equation}

In view of the above observation, now we give the crucial identity as follows

\begin{proposition}\label{Pro-1}
\begin{equation}\label{2.7}
\begin{split}
\,& \mathbf{Re}Z_{\overline{\alpha}}\Big{\{}e^{2(n-1)f}\Big[\big(g+3\sqrt{-1}f_0\big)E_{\alpha}\\
\,&\hspace{77pt}      +\big(g- \sqrt{-1}f_0\big)D_{\alpha}     -3\sqrt{-1}f_0 G_{\alpha}
        -\frac{p}{4}f_{\alpha}|\partial f|^4 \Big]\Big{\}}\\
 =&\, e^{(2n+p)f}\big(|E_{\alpha\overline{\beta}}|^2  +|D_{\alpha\beta}|^2\big)\\
 \,&  +e^{2(n-1)f}\big(|G_{\alpha}|^2+|G_{\alpha}+D_{\alpha}|^2  +|G_{\alpha}-E_{\alpha}|^2
                       +|D_{\alpha\beta}f_{\overline{\gamma}}+E_{\alpha\overline{\gamma}}f_{\beta}|^2\big)\\
 \, & +e^{(2n-2)f}\mathbf{Re}\big(D_{\alpha}+E_{\alpha}\big)  f_{\overline{\alpha}}\big( pe^{(2+p)f} -\frac{p}{2}|\partial f|^2 \big)\\
 \,& -p(2n-1)|\partial f|^2e^{2(n+1+p)f}-\frac{p}{4}(7n-6)|\partial f|^4e^{(2n+p)f}\\
 \,& -\frac{p}{4}n|\partial f|^6e^{2(n-1)f} -3np|f_0|^2e^{(2n+p)f}.
\end{split}
\end{equation}
\end{proposition}

\begin{remark}
Note that for $p=0$, then (\ref{2.7}) is exactly a remarkable identity found by Jerison and Lee (see (4.2) in  \cite{JL1988}).
\end{remark}

\vspace{10pt}
$\mathbf{Proof\,\, of\,\, proposition\,\, \ref{Pro-1}}$ \qquad Denote

$$\mathcal{L}=\mathcal{L}_1+\mathcal{L}_2+\mathcal{L}_3+\mathcal{L}_4,$$
with

$$\quad\mathcal{L}_1=Z_{\overline{\alpha}}\Big{\{} \big(g+3\sqrt{-1}f_0\big)E_{\alpha}e^{2(n-1)f} \Big{\}},$$

$$\hspace{-2pt}\quad\mathcal{L}_2=Z_{\overline{\alpha}}\Big{\{} \big(g- \sqrt{-1}f_0\big)D_{\alpha}e^{2(n-1)f} \Big{\}},$$

$$\mathcal{L}_3=Z_{\overline{\alpha}}\Big{\{} -3\sqrt{-1}f_0 G_{\alpha}e^{2(n-1)f} \Big{\}},$$

$$\hspace{-10pt} \mathcal{L}_4=Z_{\overline{\alpha}}\Big{\{} -\frac{p}{4}f_{\alpha}|\partial f|^4e^{2(n-1)f} \Big{\}}.$$

First we compute $\mathcal{L}_3$. We have, by (\ref{2.4}) and  the commutative formulae,

\begin{equation}\label{2.8}
\begin{split}
 G_{\alpha,\overline{\alpha}}
 =&\, \sqrt{-1}f_{0\alpha\overline{\alpha}}+g_{\overline{\alpha}}f_{\alpha} +gf_{\alpha\overline{\alpha}}\\
 =&\, \sqrt{-1}f_{\alpha\overline{\alpha} 0}+g_{\overline{\alpha}}f_{\alpha} +g(f_{\overline{\alpha}\alpha}+2n\sqrt{-1}f_0)\\
 =&\, f_{\alpha}g_{\overline{\alpha}}-n\sqrt{-1}g_{0}-n|g|^2+2n\sqrt{-1}f_0g.
\end{split}
\end{equation}
There for

\begin{equation}\label{2.9}
\begin{split}
 e^{-2(n-1)f}\mathcal{L}_3
 =&\, e^{-2(n-1)f}Z_{\overline{\alpha}}\Big{\{} -3\sqrt{-1}f_0 G_{\alpha}e^{2(n-1)f} \Big{\}}\\
 =&\, -3\sqrt{-1}f_{0}G_{\alpha,\overline{\alpha}}-3\sqrt{-1}f_{0\overline{\alpha}}G_{\alpha}
                                                 -6(n-1)\sqrt{-1}f_{0}f_{\overline{\alpha}}G_{\alpha}\\
 =&\, -3\sqrt{-1}f_{0}\big( f_{\alpha}g_{\overline{\alpha}}-n\sqrt{-1}g_{0}-n|g|^2+2n\sqrt{-1}f_0g \big)\\
\,&\, +3\sqrt{-1}(G_{\overline{\alpha}}-\overline{g}f_{\overline{\alpha}})G_{\alpha} -6(n-1)\sqrt{-1}f_{0}f_{\overline{\alpha}}G_{\alpha}\\
 =&\, 3|G_{\alpha}|^2 -3(\overline{g}+2(n-1)\sqrt{-1}f_0)f_{\overline{\alpha}}G_{\alpha}\\
\,&\, -3\sqrt{-1}f_0f_{\alpha}g_{\overline{\alpha}}-3nf_0g_{0}+3n\sqrt{-1}f_0|g|^2+6n|f_0|^2g.
\end{split}
\end{equation}

Next we compute $\mathcal{L}_1$. Also by (\ref{2.4}) and  the commutative formulae,

\begin{equation}\label{2.10}
\begin{split}
 E_{\alpha,\overline{\alpha}}
 =&\, f_{\alpha\overline{\beta}\overline{\alpha}}f_{\beta}+f_{\alpha\overline{\beta}}f_{\beta\overline{\alpha}}
      +g_{\overline{\alpha}}f_{\alpha} +gf_{\alpha\overline{\alpha}}\\
 =&\, f_{\alpha\overline{\alpha}\overline{\beta}}f_{\beta}
       +f_{\alpha\overline{\beta}}(f_{\overline{\alpha}\beta}+2\sqrt{-1}f_0\delta_{\beta\overline{\alpha}})
       +g_{\overline{\alpha}}f_{\alpha} +gf_{\alpha\overline{\alpha}}\\
 =&\, -ng_{\overline{\beta}}f_{\beta} +f_{\alpha\overline{\beta}}f_{\overline{\alpha}\beta}
       +2\sqrt{-1}f_0f_{\alpha\overline{\alpha}}
       +g_{\overline{\alpha}}f_{\alpha} +gf_{\alpha\overline{\alpha}}\\
 =&\, (1-n)f_{\alpha}g_{\overline{\alpha}}
      +(E_{\alpha\overline{\beta}}-g\delta_{\alpha\overline{\beta}})(E_{\overline{\alpha}\beta}-\overline{g}\delta_{\overline{\alpha}\beta})
      -n|g|^2\\
 =&\, |E_{\alpha\overline{\beta}}|^2+(1-n)f_{\alpha}g_{\overline{\alpha}}.
\end{split}
\end{equation}
There for

\begin{equation}\label{2.11}
\begin{split}
 e^{-2(n-1)f}\mathcal{L}_1
 =&\, e^{-2(n-1)f}Z_{\overline{\alpha}}\Big{\{} \big(g+3\sqrt{-1}f_0\big)E_{\alpha}e^{2(n-1)f} \Big{\}}\\
 =&\, \big(g+3\sqrt{-1}f_0\big)E_{\alpha,\overline{\alpha}}\\
\,&\, +\big(g_{\overline{\alpha}}+3\sqrt{-1}f_{0\overline{\alpha}}\big)E_{\alpha}
      +2(n-1)\big(g+3\sqrt{-1}f_0\big)f_{\overline{\alpha}}E_{\alpha}\\
 =&\, \big(g+3\sqrt{-1}f_0\big)\big( |E_{\alpha\overline{\beta}}|^2+(1-n)f_{\alpha}g_{\overline{\alpha}} \big)\\
\,&\, +g_{\overline{\alpha}}E_{\alpha}+3\big(-G_{\overline{\alpha}}+\overline{g}f_{\overline{\alpha}}\big)E_{\alpha}
      +2(n-1)\big(g+3\sqrt{-1}f_0\big)f_{\overline{\alpha}}E_{\alpha}\\
 =&\, \big(g+3\sqrt{-1}f_0\big)|E_{\alpha\overline{\beta}}|^2+\big( g_{\overline{\alpha}}-3G_{\overline{\alpha}}\big)E_{\alpha}\\
\,&\, +\big(3\overline{g}+2(n-1)(g+3\sqrt{-1}f_0)\big)f_{\overline{\alpha}}E_{\alpha}\\
\,&\, +(1-n)(g+3\sqrt{-1}f_0)f_{\alpha}g_{\overline{\alpha}}.
\end{split}
\end{equation}

Now we compute $\mathcal{L}_2$. Using the commutative formulae, we compute

\begin{equation}\label{2.12}
\begin{split}
 f_{\alpha\beta\overline{\alpha}}
 =&\, f_{\alpha\overline{\alpha}\beta}+2\sqrt{-1}f_{0\alpha}\delta_{\beta\overline{\alpha}}\\
 =&\, (f_{\overline{\alpha}\alpha}+2n\sqrt{-1}f_0)_{,\beta}+2\sqrt{-1}f_{0\beta}\\
 =&\, -n\overline{g}_{\beta}+2(n+1)(G_{\beta}-gf_{\beta})\\
 =&\, 2(n+1)G_{\beta}-n\overline{g}_{\beta}-2(n+1)f_{\beta}g.
\end{split}
\end{equation}
By this and (\ref{2.4}), (\ref{2.5}) we deduce

\begin{equation}\label{2.13}
\begin{split}
 D_{\alpha,\overline{\alpha}}
 =&\, f_{\alpha\beta\overline{\alpha}}f_{\overline{\beta}}+f_{\alpha\beta}f_{\overline{\beta}\overline{\alpha}}
      -2(|\partial f|^2)_{,\overline{\alpha}}f_{\alpha} -2|\partial f|^2f_{\alpha\overline{\alpha}}\\
 =&\, \big( 2(n+1)G_{\beta}-n\overline{g}_{\beta}-2(n+1)f_{\beta}g \big)f_{\overline{\beta}}\\
\,&\, +(D_{\alpha\beta}+2f_{\alpha}f_{\beta})(D_{\overline{\alpha}\overline{\beta}}+2f_{\overline{\alpha}}f_{\overline{\beta}})\\
\,&\, -2\big( D_{\overline{\alpha}}+E_{\overline{\alpha}}+\overline{g}f_{\overline{\alpha}}-2f_{\overline{\alpha}}e^{(2+p)f} \big)f_{\alpha}
      +2n|\partial f|^2g\\
 =&\, |D_{\alpha\beta}|^2 +2f_{\overline{\alpha}}D_{\alpha} -2f_{\alpha}E_{\overline{\alpha}}
      +2(n+1)f_{\overline{\alpha}}G_{\alpha} -nf_{\overline{\alpha}}\overline{g}_{\alpha}.
\end{split}
\end{equation}
There for

\begin{equation}\label{2.14}
\begin{split}
\,&   e^{-2(n-1)f}\mathcal{L}_2\\
 =&\,\, e^{-2(n-1)f}Z_{\overline{\alpha}}\Big{\{} \big(g- \sqrt{-1}f_0\big)D_{\alpha}e^{2(n-1)f} \Big{\}}\\
 =&\,\, \big(g- \sqrt{-1}f_0\big)D_{\alpha,\overline{\alpha}}\\
\,&\,\, +\big(g_{\overline{\alpha}}- \sqrt{-1}f_{0\overline{\alpha}}\big)D_{\alpha}
        +2(n-1)\big(g- \sqrt{-1}f_0\big)f_{\overline{\alpha}}D_{\alpha}\\
 =&\,\, \big(g- \sqrt{-1}f_0\big)\big( |D_{\alpha\beta}|^2 +2f_{\overline{\alpha}}D_{\alpha} -2f_{\alpha}E_{\overline{\alpha}}
                                       +2(n+1)f_{\overline{\alpha}}G_{\alpha} -nf_{\overline{\alpha}}\overline{g}_{\alpha} \big)\\
\,&\,\, +\big(g_{\overline{\alpha}}+G_{\overline{\alpha}}-\overline{g}f_{\overline{\alpha}} \big)D_{\alpha}
        +2(n-1)\big(g- \sqrt{-1}f_0\big)f_{\overline{\alpha}}D_{\alpha}\\
 =&\,\, \big(g- \sqrt{-1}f_0\big)|D_{\alpha\beta}|^2+(g_{\overline{\alpha}}+G_{\overline{\alpha}})D_{\alpha}\\
\,&\,\, +(2ng-\overline{g}-2n\sqrt{-1}f_0)f_{\overline{\alpha}}D_{\alpha} -2(g-\sqrt{-1}f_0)f_{\alpha}E_{\overline{\alpha}}\\
\,&\,\, +2(n+1)(g-\sqrt{-1}f_0)f_{\overline{\alpha}}G_{\alpha} -n(g-\sqrt{-1}f_0)f_{\overline{\alpha}}\overline{g}_{\alpha} .
\end{split}
\end{equation}

Finally, for $\mathcal{L}_4 $, by (\ref{2.3}) and (\ref{2.5}), a direct computation shows

\begin{equation}\label{2.15}
\begin{split}
   e^{-2(n-1)f}\mathcal{L}_4
 =&     e^{-2(n-1)f}Z_{\overline{\alpha}}\Big{\{} -\frac{p}{4}f_{\alpha}|\partial f|^4e^{2(n-1)f} \Big{\}}\\
 =&     -\frac{p}{2}\big( E_{\overline{\alpha}}+D_{\overline{\alpha}}\big)f_{\alpha}|\partial f|^2\\
\,&\,\, -\frac{p}{4}n|\partial f|^6  +\frac{p}{4}(n+2)|\partial f|^4e^{(2+p)f} -\frac{p}{4}(n+2)\sqrt{-1}f_0|\partial f|^4.
\end{split}
\end{equation}

By (\ref{2.9}), (\ref{2.11}), (\ref{2.14}) and (\ref{2.15}), noticing
$f_{\alpha}E_{\overline{\alpha}}=f_{\overline{\alpha}}E_{\alpha} $ (this also implies it's real), we obtain

\begin{equation}\label{2.16}
\begin{split}
\,&   e^{-2(n-1)f}(\mathcal{L}_1+\mathcal{L}_2+\mathcal{L}_3+\mathcal{L}_4)\\
 =&\,\, \big(g- \sqrt{-1}f_0\big)|D_{\alpha\beta}|^2 +\big(g+3\sqrt{-1}f_0\big)|E_{\alpha\overline{\beta}}|^2 +3|G_{\alpha}|^2\\
\,&\,\, +(g_{\overline{\alpha}}+G_{\overline{\alpha}})D_{\alpha}\,
        +\big( g_{\overline{\alpha}}-3G_{\overline{\alpha}}\big)E_{\alpha}\\
\,&\,\, +\Big( 2ng-\overline{g}-\frac{p}{2}|\partial f|^2-2n\sqrt{-1}f_0 \Big)f_{\overline{\alpha}}D_{\alpha}\\
\,&\,\, +\Big( 2(n-2)g+3\overline{g}-\frac{p}{2}|\partial f|^2+(6n-4)\sqrt{-1}f_0) \Big)f_{\overline{\alpha}}E_{\alpha}\\
\,&\,\, +\Big( 2(n+1)g-3\overline{g}-(8n-4)\sqrt{-1}f_0 \Big)f_{\overline{\alpha}}G_{\alpha}\\
\,&\,\, -\big( (n-1)g+3n\sqrt{-1}f_0 \big)f_{\alpha}g_{\overline{\alpha}}\,
        -n(g-\sqrt{-1}f_0)f_{\overline{\alpha}}\overline{g}_{\alpha} -3nf_0g_0\\
\,&\,\, +3n\sqrt{-1}f_0|g|^2 +6n|f_0|^2g \\
\,&\,\, -\frac{p}{4}n|\partial f|^6  +\frac{p}{4}(n+2)|\partial f|^4e^{(2+p)f} -\frac{p}{4}(n+2)\sqrt{-1}f_0|\partial f|^4.
\end{split}
\end{equation}
Straight calculations show

\begin{equation}\label{2.17}
\begin{split}
g_0
 =&\,\, \sqrt{-1}f_{\alpha}G_{\overline{\alpha}}-\sqrt{-1}f_{\overline{\alpha}}G_{\alpha}\\
\,&\,\, +2f_0|\partial f|^2 +(2+p)f_0e^{(2+p)f} -\sqrt{-1}f_{00}.
\end{split}
\end{equation}
By this and virtue of (\ref{2.6}), we finally reach

\begin{equation}\label{2.18}
\begin{split}
\,&   e^{-2(n-1)f}(\mathcal{L}_1+\mathcal{L}_2+\mathcal{L}_3+\mathcal{L}_4)\\
 =&\,\, \big(g- \sqrt{-1}f_0\big)|D_{\alpha\beta}|^2 +\big(g+3\sqrt{-1}f_0\big)|E_{\alpha\overline{\beta}}|^2 +3|G_{\alpha}|^2\\
\,&\,\, +(D_{\overline{\alpha}}+E_{\overline{\alpha}}+2G_{\overline{\alpha}})D_{\alpha}
        +\big( D_{\overline{\alpha}}+E_{\overline{\alpha}}-2G_{\overline{\alpha}}\big)E_{\alpha}\\
\,&\,\, +p\big( f_{\overline{\alpha}} -\frac{1}{2}f_{\overline{\alpha}}|\partial f|^2 \big)\big(D_{\alpha}+E_{\alpha}\big)\\
\,&\,\, +(n-1)\big( |\partial f|^2+e^{(2+p)f} \big)(f_{\overline{\alpha}}G_{\alpha}-f_{\alpha}G_{\overline{\alpha}})\\
\,&\,\, +(5n+1)\sqrt{-1}f_0(f_{\overline{\alpha}}G_{\alpha}+f_{\alpha}G_{\overline{\alpha}})\\
\,&\,\, -p(2n-1)|\partial f|^2e^{2(2+p)f}-\frac{p}{4}(7n-6)|\partial f|^4e^{(2+p)f}\\
\,&\,\, -\frac{p}{4}n|\partial f|^6 -3np|f_0|^2e^{(2+p)f}\\
\,&\,\, -\frac{p}{4}(n+2)\sqrt{-1}f_0|\partial f|^4-p\sqrt{-1}f_0|\partial f|^2e^{(2+p)f}\\
\,&\,\, +3n\sqrt{-1}f_0|g|^2 -6n\sqrt{-1}f_0|f_0|^2+3n\sqrt{-1}f_0f_{00}.
\end{split}
\end{equation}
From this one can check (\ref{2.7}) easily and complete the proof of proposition \ref{Pro-1}. \qed

\section{Proof of Theorem \ref{Thm3} }

\setcounter{equation}{0}
\setcounter{theorem}{0}

$\mathbf{Proof\,\, of\,\, Theorem\,\, \ref{Thm3}}$\qquad Let $f$ satisfy the equation (\ref{2.3}) and hence the identity (\ref{2.7}).
  Then by $q=q^*+ \frac{p}{n}$, the subcritical exponent $1<q<q^*$ corresponding to $-2<p< 0$ and (\ref{1.3}) is equivalent to

\begin{equation}\label{3-1}
 \int_{B_r(\xi_0)} e^{(2n+4+3p)f}
 \lesssim \,\,C r^{2n+2-2\times\frac{2n+4+3p}{2+p}}.
\end{equation}
 Note that (\ref{2.7})  can be rewritten as

\begin{equation}\label{3-2}
\begin{split}
\mathcal{M}= \mathbf{Re}Z_{\overline{\alpha}}\Big{\{} \Big[ &\big(D_{\alpha}+E_{\alpha})(|\partial f|^2+e^{(2+p)f}) \\
                           \,&-\sqrt{-1}f_0\big(2D_{\alpha}-2E_{\alpha}+3 G_{\alpha}\big)
                              -\frac{p}{4}f_{\alpha}|\partial f|^4 \Big]e^{2(n-1)f}\Big{\}},
\end{split}
\end{equation}
with
\begin{equation}\label{3-3}
\begin{split}
\mathcal{M}
  =& \big(|E_{\alpha\overline{\beta}}|^2  +|D_{\alpha\beta}|^2\big)e^{(2n+p)f}
     +\Big( |G_{\alpha}|^2+|D_{\alpha\beta}f_{\overline{\gamma}}+E_{\alpha\overline{\gamma}}f_{\beta}|^2\Big)e^{2(n-1)f}\\
 \,&  + s\Big(|G_{\alpha}+D_{\alpha}|^2 + |G_{\alpha}-E_{\alpha}|^2\Big)e^{2(n-1)f}\\
 \,&  +(1-s)\Big(|G_{\alpha}+D_{\alpha}|^2 + |G_{\alpha}-E_{\alpha}|^2\Big)e^{2(n-1)f}\\
 \,& + p e^{2(n-1)f}\mathbf{Re}[f_{\overline{\alpha}}(G_{\alpha}+D_{\alpha})]\big( e^{(2+p)f} -\frac{1}{2}|\partial f|^2 \big)\\
 \,& + p e^{2(n-1)f}\mathbf{Re}[f_{\overline{\alpha}}(E_{\alpha}-G_{\alpha})]\big( e^{(2+p)f} -\frac{1}{2}|\partial f|^2 \big)\\
 \,& -p(2n-1)|\partial f|^2e^{2(n+1+p)f}-\frac{p}{4}(7n-6)|\partial f|^4e^{(2n+p)f}\\
 \,& -\frac{p}{4}n|\partial f|^6e^{2(n-1)f} -3np|f_0|^2e^{(2n+p)f}.
\end{split}
\end{equation}
So for $-2<p < 0$, we shall choose suitable $0<s<1$ such that  $\mathcal{M}\geq 0$.

Now we rewrite $\mathcal{M}$ as
\begin{equation}\label{3-3a}
\begin{split}
\mathcal{M}
  =& \big(|E_{\alpha\overline{\beta}}|^2  +|D_{\alpha\beta}|^2\big)e^{(2n+p)f}
     +\Big( |G_{\alpha}|^2+|D_{\alpha\beta}f_{\overline{\gamma}}+E_{\alpha\overline{\gamma}}f_{\beta}|^2\Big)e^{2(n-1)f}\\
 \,& + s \Big(|G_{\alpha}+D_{\alpha}|^2 + |G_{\alpha}-E_{\alpha}|^2\Big)e^{2(n-1)f}\\
 \,& +e^{2(n-1)f}\Big|\sqrt{1-s}(G_{\alpha}+D_{\alpha}) +\frac{p}{2\sqrt{1-s}} f_{\alpha} \big( e^{(2+p)f} -\frac{1}{2}|\partial f|^2 \big)\Big|^2\\
 \,& +e^{2(n-1)f}\Big|\sqrt{1-s}(E_{\alpha}-G_{\alpha}) +\frac{p}{2\sqrt{1-s}} f_{\alpha} \big( e^{(2+p)f} -\frac{1}{2}|\partial f|^2 \big)\Big|^2\\
 \,&- \frac{p^2}{2(1-s)}e^{2(n-1)f}|\partial f|^2\Big[ e^{2(2+p)f} + \frac{1}{4}|\partial f|^4 -  e^{(2+p)f}|\partial f|^2 \Big]\\
 \,& -p(2n-1)|\partial f|^2e^{2(n+1+p)f}-\frac{p}{4}(7n-6)|\partial f|^4e^{(2n+p)f}\\
 \,& -\frac{p}{4}n|\partial f|^6e^{2(n-1)f} -3np|f_0|^2e^{(2n+p)f}.
\end{split}
\end{equation}

 Then we treat the terms in the last three lines and  get
\begin{equation}\label{3-3b}
\begin{split}
\mathcal{M}
  =& \big(|E_{\alpha\overline{\beta}}|^2  +|D_{\alpha\beta}|^2\big)e^{(2n+p)f}
     +\Big( |G_{\alpha}|^2+|D_{\alpha\beta}f_{\overline{\gamma}}+E_{\alpha\overline{\gamma}}f_{\beta}|^2\Big)e^{2(n-1)f}\\
    \,& + s \Big(|G_{\alpha}+D_{\alpha}|^2 + |G_{\alpha}-E_{\alpha}|^2\Big)e^{2(n-1)f}\\
   \,& +e^{2(n-1)f}\Big|\sqrt{1-s} (G_{\alpha}+D_{\alpha}) + \frac{p}{2\sqrt{1-s}} f_{\alpha} \big( e^{(2+p)f} -\frac{1}{2}|\partial f|^2 \big)\Big|^2\\
   \,& +e^{2(n-1)f}\Big|\sqrt{1-s} (E_{\alpha}-G_{\alpha}) + \frac{p}{2\sqrt{1-s}} f_{\alpha} \big( e^{(2+p)f} -\frac{1}{2}|\partial f|^2 \big)\Big|^2\\
    \,& -p\big[\frac{n}{4} + \frac{p}{8(1-s)}\big]|\partial f|^6e^{2(n-1)f}-\frac{p}{4}\big[7n-6 -\frac{2p}{1-s}\big] |\partial f|^4e^{(2n+p)f}\\
    \,& -p\big[2n-1+\frac{p}{2(1-s)}\big]|\partial f|^2e^{2(n+1+p)f} -3np|f_0|^2e^{(2n+p)f}.
\end{split}
\end{equation}

Now we take $0<s=s_0= \frac{1}{2} +\frac{p}{4n}<1$, then
\begin{equation}\label{3-3c}
\begin{split}
\mathcal{M}
  =& \big(|E_{\alpha\overline{\beta}}|^2  +|D_{\alpha\beta}|^2\big)e^{(2n+p)f}
     +\Big( |G_{\alpha}|^2+|D_{\alpha\beta}f_{\overline{\gamma}}+E_{\alpha\overline{\gamma}}f_{\beta}|^2\Big)e^{2(n-1)f}\\
 \,& + s_0 \Big(|G_{\alpha}+D_{\alpha}|^2 + |G_{\alpha}-E_{\alpha}|^2\Big)e^{2(n-1)f}\\
 \,& +e^{2(n-1)f}\Big|\sqrt{1-s_0}(G_{\alpha}+D_{\alpha}) +\frac{p}{2\sqrt{1-s_0}} f_{\alpha}\big( e^{(2+p)f} -\frac{1}{2}|\partial f|^2 \big)\Big|^2\\
 \,& +e^{2(n-1)f}\Big|\sqrt{1-s_0}(E_{\alpha}-G_{\alpha}) +\frac{p}{2\sqrt{1-s_0}} f_{\alpha}\big( e^{(2+p)f} -\frac{1}{2}|\partial f|^2 \big)\Big|^2\\
 \,& -p\frac{n(2n+p)}{4(2n-p)}|\partial f|^6e^{2(n-1)f}-\frac{p}{4}\big[7n-6 -\frac{8np}{2n-p}\big]  |\partial f|^4e^{(2n+p)f}\\
 \,& -p\frac{4n^2-2n+p}{2n-p}|\partial f|^2e^{2(n+1+p)f} -3np|f_0|^2e^{(2n+p)f},
\end{split}
\end{equation}
and clearly all the coefficients in above are positive for $-2<p<0$ and  $\mathcal{M}\geq 0$.

Since $B_{4r}\subset \Omega$, we can take  a real smooth cut off function  $\eta$  such that

\begin{equation}\label{3-4}
\begin{cases}
                          \eta\equiv 1  &in \,\,B_r,\\
                        0\leq\eta\leq1  &in \,\,B_{2r},\\
                          \eta\equiv 0  &in \,\,\Omega\backslash B_{2r},\\
     |\partial \eta|\lesssim \frac{1}{r}  &in \,\,\Omega,
\end{cases}
\end{equation}
\noindent where we use ``$\lesssim $" , ``$\cong$"  to replace ``$\leq$" and ``$=$" respectively, to drop out some
positive constants independent of $r$ and $f$.

Take a real $s>0$ big enough. Multiply both sides of (\ref{3-2}) by $\eta^s $ and integrate over  $\Omega$ we have

\begin{equation}\label{3-5}
\begin{split}
\,& \int_{\Omega}\eta^s \mathcal{M}\\
 =& \int_{\Omega}\eta^s\mathbf{Re}Z_{\overline{\alpha}}\Big{\{} \big[ \big(D_{\alpha}+E_{\alpha})(|\partial f|^2+e^{(2+p)f}) \\
\,&\hspace{68pt} -\sqrt{-1}f_0\big(2D_{\alpha}-2E_{\alpha}+3 G_{\alpha}\big)
                 -\frac{p}{4}f_{\alpha}|\partial f|^4 \big]e^{2(n-1)f}\Big{\}}.
\end{split}
\end{equation}
Integrating by part and using (\ref{3-4}) we get

\begin{equation}\label{3-6}
\begin{split}
\,& \int_{\Omega}\eta^s \mathcal{M}\\
= & -s\int_{\Omega}\eta^{s-1}\mathbf{Re}\eta_{\overline{\alpha}}\Big{\{} \big[ \big(D_{\alpha}+E_{\alpha})(|\partial f|^2+e^{(2+p)f}) \\
\,&\hspace{99pt} -\sqrt{-1}f_0\big(2D_{\alpha}-2E_{\alpha}+3 G_{\alpha}\big)
                 -\frac{p}{4}f_{\alpha}|\partial f|^4 \big]e^{2(n-1)f}\Big{\}}\\
\lesssim & \frac{1}{r}\int_{\Omega}\eta^{s-1}\Big{\{} |D_{\alpha}+E_{\alpha}|(|\partial f|^2+e^{(2+p)f})e^{2(n-1)f} \\
\,& \hspace{55pt}    +|f_0|\big|2D_{\alpha}-2E_{\alpha}+3 G_{\alpha}\big|e^{2(n-1)f}+|\partial f|^5e^{2(n-1)f}\Big{\}}
\end{split}
\end{equation}
Since

$$|D_{\alpha}+E_{\alpha}|\leq |D_{\alpha}+G_{\alpha}|+|E_{\alpha}-G_{\alpha}|,$$

$$\big|2D_{\alpha}-2E_{\alpha}+3 G_{\alpha}\big| \leq 2|D_{\alpha}+G_{\alpha}|+2|E_{\alpha}-G_{\alpha}|+|G_{\alpha}|,$$

\noindent using the Young's inequality $ab\leq \epsilon a^2+\frac{C}{\epsilon}b^2$ in (\ref{3-6}) we obtain

\begin{equation}\label{3-7}
\begin{split}
\int_{\Omega}\eta^s \mathcal{M}
\lesssim & \epsilon \int_{\Omega}\eta^{s}\big( |D_{\alpha}+G_{\alpha}|^2+|E_{\alpha}-G_{\alpha}|^2+|G_{\alpha}|^2\big) e^{2(n-1)f} \\
\,& \hspace{13pt}  +\frac{1}{\epsilon r^2} \int_{\Omega}\eta^{s-2}\big(|\partial f|^4+e^{2(2+p)f}+|f_0|^2\big)e^{2(n-1)f}\\
\,& \hspace{113pt}  +\frac{1}{r} \int_{\Omega}\eta^{s-1}|\partial f|^5e^{2(n-1)f}.
\end{split}
\end{equation}
This implies, by taking $\epsilon$ small, that

\begin{equation}\label{3-8}
\begin{split}
\int_{\Omega}\eta^s \mathcal{M}
\lesssim & \frac{1}{r^2} \int_{\Omega}\eta^{s-2}\big(|\partial f|^4+e^{2(2+p)f}+|f_0|^2\big)e^{2(n-1)f}\\
\,&        \hspace{86pt}  +\frac{1}{r} \int_{\Omega}\eta^{s-1}|\partial f|^5e^{2(n-1)f}.
\end{split}
\end{equation}

 To go forward, we need the following  lemmas, that will be proved at the end of this section.

\begin{lemma}\label{lem-1}
\begin{equation}\label{3-9}
\begin{split}
 \int_{\Omega}\eta^{s-2}|f_0|^2e^{2(n-1)f}
\lesssim &\, \epsilon r^2\int_{\Omega} \eta^s \mathcal{M}
            +\int_{\Omega}\eta^{s-2}|\partial f|^4e^{2(n-1)f}\\
\,&         +\int_{\Omega}\eta^{s-2}|\partial f|^2e^{(2n+p)f}
            +\frac{1}{r^2}\int_{\Omega}\eta^{s-4}|\partial f|^2e^{2(n-1)f}.
\end{split}
\end{equation}
\end{lemma}

\begin{lemma}\label{lem-2}
\begin{equation}\label{3-10}
 \int_{\Omega} \eta^s e^{(2n+4+3p)f}
 \lesssim \, \int_{\Omega} \eta^s|\partial f|^2e^{2(n+1+p)f}
       +\frac{1}{r^{2}}\int_{\Omega} \eta^{s-2} e^{(2n+2+2p)f} .
\end{equation}
\end{lemma}

\vspace{10pt}
 Now plugging (\ref{3-9}) into  (\ref{3-8}) with small $\epsilon$ we get

\begin{equation}\label{3-11}
\begin{split}
\int_{\Omega}\eta^s \mathcal{M}
\lesssim &\,\frac{1}{r^{2}}\int_{\Omega}\eta^{s-2} e^{2(n+1+p)f}\\
\,&        +\frac{1}{r^2}\int_{\Omega}\eta^{s-2}|\partial f|^4e^{2(n-1)f}
           +\frac{1}{r^2}\int_{\Omega}\eta^{s-2}|\partial f|^2e^{(2n+p)f}\\
\,&        +\frac{1}{r^4}\int_{\Omega}\eta^{s-4}|\partial f|^2e^{2(n-1)f}
           +\frac{1}{r} \int_{\Omega}\eta^{s-1}|\partial f|^5e^{2(n-1)f}.
\end{split}
\end{equation}
For the last term in above, using Young's inequality one get

\begin{equation}\label{3-12}
\begin{split}
\frac{1}{r} \int_{\Omega}\eta^{s-1}|\partial f|^5e^{2(n-1)f}
\lesssim &    \epsilon \int_{\Omega}\eta^{s }|\partial f|^6e^{2(n-1)f}
             +\frac{1}{r^6}\int_{\Omega}\eta^{s-6} e^{2(n-1)f}.
\end{split}
\end{equation}
Similarly, one has

\begin{equation}\label{3-13}
\begin{split}
\frac{1}{r^2} \int_{\Omega}\eta^{s-2}|\partial f|^4e^{2(n-1)f}
\lesssim &  \epsilon \int_{\Omega}\eta^{s }|\partial f|^6e^{2(n-1)f}
            +\frac{1}{r^6}\int_{\Omega}\eta^{s-6} e^{2(n-1)f},
\end{split}
\end{equation}

\begin{equation}\label{3-14}
\begin{split}
\frac{1}{r^2} \int_{\Omega}\eta^{s-2}|\partial f|^2e^{(2n+p)f}
\lesssim & \epsilon \int_{\Omega}\eta^{s }|\partial f|^4e^{(2n+p)f}
             +\frac{1}{r^4}\int_{\Omega}\eta^{s-4} e^{(2n+p)f},
\end{split}
\end{equation}
and
\begin{equation}\label{3-15}
\begin{split}
\frac{1}{r^4} \int_{\Omega}\eta^{s-4}|\partial f|^2e^{2(n-1)f}
\lesssim &   \epsilon \int_{\Omega}\eta^{s }|\partial f|^6e^{2(n-1)f}
             +\frac{1}{r^6}\int_{\Omega}\eta^{s-6} e^{2(n-1)f}.
\end{split}
\end{equation}
Inserting these into (\ref{3-11}) and taking $\epsilon$ small yields

\begin{equation}\label{3-16}
\begin{split}
\int_{\Omega}\eta^s \mathcal{M}
\lesssim &  \frac{1}{r^{2}}\int_{\Omega}\eta^{s-2} e^{2(n+1+p)f}\\
\, &        +\frac{1}{r^4}\int_{\Omega}\eta^{s-4} e^{(2n+p)f}+ \frac{1}{r^{6}}\int_{\Omega}\eta^{s-6} e^{2(n-1)f}.
\end{split}
\end{equation}
Combining this with (\ref{3-10}) we arrive at

\begin{equation}\label{3-17}
\begin{split}
\,& \int_{\Omega} \eta^s e^{(2n+4+3p)f}\\
 \lesssim &  \frac{1}{r^{2}}\int_{\Omega}\eta^{s-2} e^{2(n+1+p)f}\\
\, &        +\frac{1}{r^4}\int_{\Omega}\eta^{s-4} e^{(2n+p)f}+ \frac{1}{r^{6}}\int_{\Omega}\eta^{s-6} e^{2(n-1)f} \\
 \lesssim & \epsilon \int_{\Omega} \eta^s e^{(2n+4+3p)f}
            +r^{-2\times\frac{2n+4+3p}{2+p}}\int_{\Omega} \eta^{s-2\times\frac{2n+4+3p}{2+p}} ,
\end{split}
\end{equation}
where in the last step, the Young's inequality has been  used three time with different exponent pairs.
Note that $0\leq\eta\leq 1$ in $\Omega$ and  $\eta= 1$ in $B_r(\xi_0)\subset\Omega$.
 Therefor, by choosing  $s>0$ big enough and $\epsilon$ small, we finally obtain

\begin{equation}\label{3-18}
 \int_{B_r(\xi_0)} e^{(2n+4+3p)f}
 \lesssim \,\, r^{2n+2-2\times\frac{2n+4+3p}{2+p}}  .
\end{equation}
This is (\ref{3-1}), and hence  theorem \ref{Thm3} is proved.\qed

\vspace{20pt}

 To complete this section, now we give the proofs of lemma \ref{lem-1},\ref{lem-2}.

\vspace{10pt}
 $\mathbf{Proof\,\, of\,\, lemma\,\, \ref{lem-1} }$

 Since $f$ satisfies the equation (\ref{2.3}), a straight calculation shows

\begin{equation}\label{3--17}
  e^{-kf}\mathbf{Re}Z_{\overline{\alpha}}\Big( \sqrt{-1}f_0f_{\alpha}e^{kf} \Big)
 = -\mathbf{Re}G_{\overline{\alpha}}f_{\alpha}-n|f_0|^2+|\partial f|^4+|\partial f|^2e^{(2+p)f}.
\end{equation}

 Multiply both sides of (\ref{3--17}) by $\eta^{s-2} e^{kf}$ with $k=2(n-1)$
 and integrate over  $\Omega$ we have

\begin{equation}\label{3.17}
\begin{split}
\,& \int_{\Omega}\eta^{s-2}\mathbf{Re}Z_{\overline{\alpha}}\Big( \sqrt{-1}f_0f_{\alpha}e^{2(n-1)f} \Big)\\
 =&\, \int_{\Omega}\eta^{s-2}\Big(-\mathbf{Re}G_{\overline{\alpha}}f_{\alpha}
                                       -n|f_0|^2+|\partial f|^4+|\partial f|^2e^{(2+p)f}\Big)e^{2(n-1)f}.
\end{split}
\end{equation}
Integrating by part, using (\ref{3-4}) and arranging the terms yields

\begin{equation}\label{3.18}
\begin{split}
 n\int_{\Omega}\eta^{s-2}|f_0|^2e^{2(n-1)f}
 =&\, \int_{\Omega}\eta^{s-2}\big(|\partial f|^4+|\partial f|^2e^{(2+p)f}\big)e^{2(n-1)f}\\
\,& \quad    -\int_{\Omega}\eta^{s-2}\mathbf{Re}G_{\overline{\alpha}}f_{\alpha}e^{2(n-1)f}\\
\,& \quad    +(s-2)\int_{\Omega}\eta^{s-3}\mathbf{Re}\eta_{\overline{\alpha}}\Big( \sqrt{-1}f_0f_{\alpha}e^{2(n-1)f} \Big)\\
\lesssim &\, \int_{\Omega}\eta^{s-2}\big(|\partial f|^4+|\partial f|^2e^{(2+p)f}\big)e^{2(n-1)f}\\
\,& \quad    +\int_{\Omega}\eta^{s-2}|G_{\overline{\alpha}}||\partial f| e^{2(n-1)f}\\
\,& \quad    +\frac{1}{r}\int_{\Omega}\eta^{s-3}|f_0||\partial f|e^{2(n-1)f}.
\end{split}
\end{equation}
For the above last two terms,  Young's inequality  implies

\begin{equation}\label{3.19}
\begin{split}
 \,&   \int_{\Omega}\eta^{s-2}|G_{\overline{\alpha}}||\partial f| e^{2(n-1)f}
       +\frac{1}{r}\int_{\Omega}\eta^{s-3}|f_0||\partial f|e^{2(n-1)f}\\
\leq &\,  \epsilon r^2\int_{\Omega}\eta^s |G_{\alpha}|^2e^{2(n-1)f}
            +\epsilon\int_{\Omega}\eta^{s-2}|f_0|^2 e^{2(n-1)f}\\
 \,& \hspace{106pt} +\frac{C}{\epsilon r^2}\int_{\Omega}\eta^{s-4}|\partial f|^2e^{2(n-1)f} .
\end{split}
\end{equation}
Submitting this into (\ref{3.18}) with small $\epsilon$  we get

\begin{equation}\label{3.20}
\begin{split}
 \int_{\Omega}\eta^{s-2}|f_0|^2e^{2(n-1)f}
\lesssim &\, \epsilon r^2\int_{\Omega} \eta^s \mathcal{M}
            +\int_{\Omega}\eta^{s-2}|\partial f|^4e^{2(n-1)f}\\
\,&         +\int_{\Omega}\eta^{s-2}|\partial f|^2e^{(2n+p)f}
            +\frac{1}{r^2}\int_{\Omega}\eta^{s-4}|\partial f|^2e^{2(n-1)f}.
\end{split}
\end{equation}
This is just (\ref{3-9}).\qed

\vspace{10pt}
 $\mathbf{Proof\,\, of\,\, lemma\,\, \ref{lem-2} }$

Multiply both sides of the equation (\ref{2.3}) by $-\eta^s e^{2(n+1+p)f}$
 and integrate over  $\Omega$ we have
\begin{equation}\label{3-21}
\begin{split}
 n\int_{\Omega} \eta^s g e^{2(n+1+p)f}
 =& -\int_{\Omega} \eta^s f_{\alpha\overline{\alpha}}e^{2(n+1+p)f}\\
 =& 2(n+1+p)\int_{\Omega} \eta^s|\partial f|^2e^{2(n+1+p)f}\\
\,& +s\int_{\Omega} \eta^{s-1} f_{\alpha}\eta_{\overline{\alpha}}e^{2(n+1+p)f}.
\end{split}
\end{equation}

Using (\ref{3-4}) and arranging the terms yields

\begin{equation}\label{3-22}
\begin{split}
 \int_{\Omega} \eta^s e^{(2n+4+3p)f}
 \lesssim & \int_{\Omega} \eta^s|\partial f|^2e^{2(n+1+p)f}+ \frac{1}{r}\int_{\Omega} \eta^{s-1} |\partial f|e^{2(n+1+p)f}\\
 \lesssim & \int_{\Omega} \eta^s|\partial f|^2e^{2(n+1+p)f}+ \frac{1}{r^2}\int_{\Omega} \eta^{s-2} e^{2(n+1+p)f},
\end{split}
\end{equation}
where in the last step, the Cauchy-Schwarz inequality has been used, and this is (\ref{3-10}) as desired. \qed

\section{Proof of theorem \ref{Thm2} }

\setcounter{equation}{0}
\setcounter{theorem}{0}

Before the proof of theorem  \ref{Thm2}, we give the following Harnack inequality,
which is a special case of that given by Capogna-Danielli-Garofalo (see Theorem 3.1 in \cite{CDG1993}),

\begin{lemma}\label{lem4.1}
Let $0\leq u\in C^2(\Omega)$  satisfies

\begin{equation}\label{4.1}
\triangle_{\mathbb{H}^{n}} u+ h(\xi) u=0 \quad \text{in}\quad \Omega.
\end{equation}
with $h(\xi)\in L^s_{loc}(\Omega)$ for some $s>\frac{Q}{2}$. Then there exist constants $C_0>0$, $r_0>0$,
such that for any $B_r(\xi)$, with $B_{4r}(\xi)\subset \Omega$, and $r<r_0$,

\begin{equation}\label{4.2}
\max_{B_r(\xi)}u\leq C_0 \min_{B_r(\xi)}u.
\end{equation}
\end{lemma}

$\mathbf{Proof\,\, of\,\, Theorem\,\, \ref{Thm2}}.$

\qquad

 Rewrite the equation (\ref{1.1}) as

\begin{equation}\label{4.3}
\triangle_{\mathbb{H}^{n}} u+ h(\xi) u=0 \quad \text{in}\quad B_1(0)\backslash\{0\},
\end{equation}
with $h(\xi)=2n^2 u^{q-1} $. For any $\xi_0\in B_1\backslash\{0\}$, take $r=\frac{1}{4}|\xi_0|$.
 Denote $ \big|B_r(\xi_0)\big|$ the volume of the ball $B_r(\xi_0)$. Using the estimate (\ref{1.3}) we have

\begin{equation}\label{4.4}
\int_{B_r(\xi_0)}h^s= (2n^2)^s\int_{B_r(\xi_0)} u^{3q-q^*} \leq C\,r^{Q-2\times\frac{3q-q^*}{q-1}}.
\end{equation}
with $s=\frac{3q-q^*}{q-1}>\frac{Q}{2}$ for $1<q<q^{*}$. This implies $h(\xi)\in L^s_{loc}(\Omega)$ for some $s>\frac{Q}{2}$
 and hence $u$ satisfies the Harnack inequality (\ref{4.2}).
 So for $|\xi_0|$ small enough, combining this Harnack inequality with (\ref{1.3}) we finally obtain

\begin{equation}\label{4.5}
\frac{1}{C}|\xi_0|^Q u(\xi_0)^{3q-q^*}\leq |B_r(\xi_0)| \Big[\frac{u(\xi_0)}{C_0}\Big]^{3q-q^*} \leq \int_{B_r(\xi_0)} u^{3q-q^*} \leq C\,r^{Q-2\times\frac{3q-q^*}{q-1}}.
\end{equation}
This implies (\ref{1.2}) and the proof of theorem \ref{Thm2} goes to the end.\qed
\begin{remark}\label{remark1}
In the following paper \cite{MaOu}, we shall generalize this methods to a class semilinear elliptic  equation on CR manifold and get a rigidity results.
\end{remark}
{\bf Acknowledgement}  Partial research of the second author was done while he was visiting
 The Chinese University of Hong Kong. He would like to thank the Institute of Mathematical Sciences
 in The Chinese University of Hong Kong for its warm hospitality.


\begin{thebibliography}{99}



\bibitem{BDN95}H. Berestycki, I. Capuzzo-Dolcetta \& L. Nirenberg, Superlinear indefinite elliptic problems and
 nonlinear Liouville theorems, Topological Methods in Nonlinear Analysis, {\bf 4} (1995), 59-78.

\bibitem{BVP01}M. F. Bidaut-Veron \& S. Pohozaev, Nonexistence results and estimates for some nonlinear elliptic problems,
J. Anal. Math. {\bf 84} (2001), 1-49.

\bibitem{BDC1997}I. Birindelli, I. Capuzzo-Dolcetta \& A. Cutri, Liouville
theorems for semilinear equations on the Heisenberg group,
Ann. Inst. Henri Poincar\'{e}-Analyse non lin\'{e}aire, {\bf 14} (1997), 295-308.

\bibitem{BP1999}I. Birindelli \& J. Prajapat, Nonlinear Liouville
theorems in the Heisenberg group via the moving plane, Comm.
Partial Differential Equations, {\bf 24} (1999), 1875-1890.

\bibitem{CDG1993}Luca Capogna, Donatella Danielli \& Nicola Garofalo,
An embedding theorem and the harnack inequality for nonlinear subelliptic equations,
Communications in Partial Differential Equations, {\bf 18:9-10} (1993), 1765-1794.

\bibitem{GV2001}N. Garofalo \& D. Vassilev, Symmetry properties of positive entire solutions of
Yamabe type equations on the groups of Heisenberg type, Duke Math. J., {\bf 106 } (2001), 411-448.

\bibitem{Ga2001}N. Gamara, CR Yamabe conjecture, the case n = 1, J. Eur. Math. Soc., {\bf 3} (2001),
105-137.

\bibitem{GY2001}N. Gamara Abdelmoula \& R. Yacoub: CR Yamabe conjecture. The case conformally flat.
Pac. J. Math., {\bf 201} (2001), No. 1, 121-175.

\bibitem{GS1981}B. Gidas \& J. Spruck, Global and local behavior of positive solutions of nonlinear
elliptic equations, Comm. Pure Appl. Math. {\bf 34} (1981), 525-598.

\bibitem{JL1987}D.S. Jerison \& J.M. Lee, The Yamabe Problem on CR manifolds, J. Differential Geom. {\bf 25} (1987), 167-197.

\bibitem{JL1988}D. Jerison \& J. M. Lee, Extremals for the Sobolev inequality on the Heisenberg group and
the CR Yamabe problem, J. Amer. Math. Soc., {\bf 1} (1988), 1-13.

\bibitem{JL1988-2}D. Jerison \& J. M. Lee, Pseudo-Einstein structures on CR manifolds, Amer. J. Math., {\bf 110} (1988), 157-178.

\bibitem{JL1989}D. Jerison \& J.M. Lee: Intrinsic CR normal coordinates and the CR Yamabe problem.
J. Diff. Geometry, {\bf 29} (1989), 303-343.

\bibitem{MaOu}X.N. Ma \& Q.Z. Ou: Rigidity results for a class semilinear elliptic on compact CR manifold, preprint.


\bibitem{Ob1971}M. Obata. The conjecture on conformal transformations of Riemannian manifolds, J. Differential
Geom., {\bf 6} (1971), 247-258.


\bibitem{Wang2015} X. Wang, On a remarkable formula of Jerison and Lee in CR geometry,
Math. Res. Lett.  {\bf 22} (2015), no. 1, 279-299.

\bibitem{Xu2009}L. Xu,  Semi-linear Liouville theorems in the Heisenberg group
via vector field methods, J. Diff. Eq., {\bf 247} (2009), 2799-2820.

\end{thebibliography}
\end{document}